\documentclass[leqno,12pt]{amsart} 
\usepackage{amssymb}
\theoremstyle{plain} 
\newtheorem{la}{Lemma}
\newtheorem*{thm}{Theorem}

\theoremstyle{definition}

\def\C{\mathbb C}
\def\N{\mathbb N}

\def\NO{{\N_0}}

\begin{document}
\title[Permutable enire functions]{Permutable
 entire functions satisfying algebraic differential equations}
\author{Walter Bergweiler}
\thanks{Supported by 
the Alexander von Humboldt Foundation
and by 
the G.I.F.,
the German--Israeli Foundation for Scientific Research and
Development, Grant G -809-234.6/2003.
}
\address{Mathematisches Seminar,
Christian--Albrechts--Universit\"at zu Kiel,
Ludewig--Meyn--Str.\ 4,
D--24098 Kiel,
Germany}
\email{bergweiler@math.uni-kiel.de}
\begin{abstract}
We show that if $f$ and $g$ are transcendental entire functions 
such that $f(g)=g(f)$, then $f$ satisfies an algebraic
differential equation if and only if $g$ does.
\end{abstract}
\maketitle

\section{Introduction}
Let $f$ and $g$ be entire or rational functions.
We say that $f$ and $g$ are {\em permutable} if $f(g)=g(f)$.
Here we prove the following result.

\begin{thm} 
Let $f$ and $g$ be permutable transcendental entire functions.
If $f$ satisfies an algebraic differential equation,
then so does $g$.
\end{thm}
We describe some related results and background.
First we note that permutability is obviously related 
to iteration since every function is clearly permutable
with all its iterates.
Fatou~\cite{Fat23} and Julia~\cite{Jul22} used the iteration 
theory they had developed to classify permutable
polynomials (and 
in fact permutable rational functions whose 
Julia sets do not coincide with  the Riemann sphere.)
The permutable rational functions were then
classified by Ritt~\cite{Rit}; a proof of Ritt's result using
iteration theory was given by Eremenko~\cite{Ere}.
The basic result obtained by these authors
says that if two rational functions 
are permutable, then they have a common iterate, except
in special cases arising from monomials, Chebychev polynomials
and the multiplication theorems of certain elliptic functions. 
These exceptional cases can be described completely.

The classification of permutable transcendental entire functions is
still an open problem.
It was shown by Bargmann~\cite[\S 4]{Bar} that
if $f$ and $g$ have a common repelling fixed point, 
then the method of Fatou and Julia can be made to work:
the Julia sets of $f$ and $g$ are equal, and if this common
Julia set does not coincide with the whole plane, then
$f$ and $g$ have a common iterate.
Note, however, that in general two permutable entire functions
need not have a common fixed point, as shown by the 
example $f(z)=z+e^z$ and $g(z)=z+e^z+2\pi i$.

For many entire functions $f$ the only nonlinear
entire functions which are permutable 
with $f$ are the iterates of $f$.
Baker~\cite{Bak58}
showed that this is the case for 
$f(z)=ae^{bz}+c$, where $a,b,c\in\C$, $a,b\neq 0$.
Ng~\cite{Ng} gave a fairly general class of entire 
functions $f$ for which any nonlinear 
entire function $g$ permutable with $f$ is of the form
$g(z)=af^n(z)+b$, where 
$a$ is a root of unity,
$b\in\C$,
$n\in\N:=\{1,2,3,\dots\}$, 
and $f^n$ is the $n$-th iterate of $f$.
Ng's class contains all functions of the form
$f(z)=e^z+p(z)$ or
$f(z)=\sin z+p(z)$ where $p$ is a nonconstant
polynomial, but it does not contain Baker's example 
$f(z)=ae^{bz}+c$.

There are several papers in which 
for certain particular functions $f$
the set of all entire 
functions $g$ of finite order which are permutable 
with $f$ is determined; see  
\cite{Kob,UY,Xio,Zhe,ZY92,ZY,ZZ}.
In some of these papers, notably~\cite{ZZ}, the fact that
the given function $f$ satisfies some differential equation
is used. In particular,
it is proved in~\cite[Theorem~1]{ZZ} that if $f$ and
$g$ are permutable entire functions of finite order,
and if $f$ satisfies a linear differential equation
with rational coefficients, 
then $g$ satisfies such a differential equation.
Our theorem above can be considered as a further contribution 
in this direction, with the restriction on the order of $g$ 
and on the linearity of the differential equation removed.

The main tools used are a result of Baker~\cite{Bak58} 
relating the growth of permutable entire functions (Lemma~\ref{baker})
and a result of Steinmetz~\cite{Ste80} concerning factorization
of solutions of certain functional and differential equations
(Lemma~\ref{steinm}).
Steinmetz's result was also used in~\cite{ZZ}.
We shall use the standard terminology of the theory of entire
and meromorphic functions as given in~\cite{Hay64}.

\section{Lemmas}
The following result due to Baker~\cite[Satz~7]{Bak58} is essential
for our argument.
\begin{la}\label{baker}
Let $f$ and $g$ be permutable transcendental entire functions.
Then there exists $p\in\N$ and $R>0$ such that 
$M(r,g)<M(r,f^p)$ for $r>R$.
\end{la}
A consequence is the following result.
\begin{la} \label{bakerT}
Let $f$ and $g$ be permutable transcendental entire functions.
Then there exists $q\in\N$ and $R>0$ such that 
$T(r,g)<T(r,f^q)$ for $r>R$.
\end{la}
{\em Proof}.
By a classical result of
P\'olya~\cite{Pol26}
(see also~\cite[Theorem~6]{Clu70} and
\cite[Theorem~2.9]{Hay64})
there
exists a constant $c>0$ such that if
$F$ and $G$ are entire, 
then
\[
M\left(r,F\left(G\right)\right)\geq M\left( c M\left(\frac{r}{2},G\right),F\right)
\]
for sufficiently large $r$. Moreover,
for entire $F$ 
we have  
\[
T(r,F) \leq \log M(r,F)\leq 3 T(2r,F)
\]
for sufficiently large $r$; 
see~\cite[Theorem~1.6]{Hay64}.
Combining these estimates we 
obtain with $p$ as in Lemma~\ref{baker}
\begin{eqnarray*}
T\left(r,f^{p+1}\right)
&\geq&
\frac13 
\log M\left(\frac{r}{2},f^{p+1}\right)\\
&\geq&
\frac13 
\log M\left( c M\left(\frac{r}{4},f\right),f^p\right)\\
&\geq&
\frac13 
\log M\left( c M\left(\frac{r}{4},f\right),g\right)\\
&\geq&
\frac13 
T\left( c M\left(\frac{r}{4},f\right),g\right)
\end{eqnarray*}
for large $r$. 
Now 
$c M(\frac{r}{4},f)>r^4$ for sufficiently large $r$
since $f$ is transcendental, and 
$T(r^4,g)\geq 3 T(r,g)$ for large $r$ since
$T(r,g)$ is convex in $\log r$.
The
conclusion follows with $q:=p+1$.\qed

Another important tool we shall use
is the following result of Steinmetz~\cite[Satz~1, Korollar~1]{Ste80}.
Generalizations and different proofs of this result have been
given by Brownawell~\cite{Bro} and Gross and Osgood
\cite{Gro89,Gro91,Gro92}.

\begin{la} \label{steinm}
Let $F_0$,$F_1$,$\dots$,$F_n$ be not identically vanishing meromorphic
functions and let
$h_0$,$h_1$,$\dots$,$h_n$ be  meromorphic functions that do not all
vanish identically. Let $g$ be a nonconstant entire function
and suppose that there exists a positive constant $K$ such that
$\sum_{j=0}^nT(r,h_j)\leq KT(r,g)$
as $r\to\infty$
outside some exceptional set of finite measure.
Suppose also that $\sum_{j=0}^n  h_j F_j(g)=0$.
Then there exist polynomials
$p_0$,$p_1$,$\dots$,$p_n$ that do not all vanish identically
such that
$\sum_{j=0}^n  p_j F_j=0$.
\end{la}

We also require some results about differential polynomials.
Let $n\in\N$,
$m_j\in \NO :=\N\cup\{0\}$
for $j=0,1,\dots,n$, and put
$m=(m_0,m_1,\dots,m_n)$.
Define $M_m[f]$ by
\[
M_m[f](z)=f(z)^{m_0}f'(z)^{m_1}f''(z)^{m_2}\dots f^{(n)}(z)^{m_n},
\]
with the convention that $M_{(0)}[f]=1$. We call
$w(m)=m_1+2m_2+\dots+nm_n$ the
{\em weight}
of $M_m[f]$.
A {\em differential polynomial}
$P[f]$ is an expression of the form
\[
\label{p}
P[f](z)=\sum_{m\in I}a_m(z)M_m[f](z),
\]
where the $a_m$ are meromorphic functions called the
{\em coefficients} of $P[f]$
and $I$ is a finite index set.
The {\em weight}
$w(P)$ of $P[f]$ is given
by $w(P)=\max_m w(m)$, where the maximum is
taken over all $m\in I$ for which $ a_m\not\equiv 0$.

Of course, 
a meromorphic function $f$ is said to satisfy an
algebraic differential equation if there exists a nontrivial
differential polynomial $P[f]$ with rational coefficients
such that $P[f]=0$.
\begin{la}\label{chain}
Let $f$ and $g$ be permutable transcendental entire functions.
Let $P[f]=\sum_{m\in I}a_m M_m[f]$ be
a differential polynomial with rational coefficients $a_m$.
Then 
$P[f](g)= \sum_{m\in I}a_m(g) M_m[f](g)=0$
 can be written in the form
$$P[f](g)= \sum_{m\in J}b_m M_m[g](f),$$ 
where the $b_m$ are meromorphic functions
satisfying $T(r,b_m)=O(U(r))$ as $r\to\infty$
outside some exceptional set of finite measure,
with 
\[
U(r)=\max\left\{ T(r,f),T(r,g)\right\}.
\]
Moreover, $\max_{m\in J}w(m)=w(P)$.
\end{la}
{\em Proof}.
Differentiating $f(g)=g(f)$ we find that
\[
f'(g)=g'(f)\frac{f'}{g'},
\]
\[
f''(g)=g''(f)\frac{(f')^2}{(g')^2}+g'(f)\frac{f''}{(g')^2}-g'(f)\frac{f'g''}{(g')^3},
\]
and so forth. 
Substituting this in $P[f](g)$ 
yields the desired representation of $P[f](g)$ with
coefficients $b_m$ which are rational
functions of $g$ and of
the derivatives of $f$ and $g$ up to order $n$.
Since $T(r,f^{(k)})\leq (1+o(1)) T(r,f)$ 
for each $k\in\N$ as $r\to\infty$ 
outside some exceptional set of finite measure (see~\cite[p.~56]{Hay64})
the conclusion follows.\qed

Finally we recall the following result
of Ostrowski~\cite[p.~269]{Ost}.
\begin{la} \label{ostr}
Let $f$ and $g$ be analytic in certain domains,
with $f$ defined in the range of $g$.
Suppose that $f$ and $g$ both satisfy some 
(possibly different) algebraic differential equation.
Then $f(g)$ satisfies some algebraic differential equation.
\end{la}
In general, of course, the differential equation for
$f(g)$ will be ``more complicated'' than those for 
$f$ and $g$; that is, its weight will be larger.

\section{Proof of the theorem}
Let $f$ and $g$ be permutable transcendental entire functions
and suppose that $f$ satisfies  an algebraic differential equation.

We choose $q$ according to  Lemma~\ref{bakerT} and put $F:=f^q$.
Lemma~\ref{ostr} implies that
$F$ also satisfies an algebraic differential equation,
say 
\[
P[F]=\sum_{m\in I}c_m M_m[F]=0,
\]
with rational functions $c_m$.
By Lemma~\ref{chain} 
we have 
$$0=P[F](g)= \sum_{m\in J}b_m M_m[g](F)$$ 
for suitable meromorphic functions
$b_m$ satisfying
$T(r,b_m)=O(U(r))$ as $r\to\infty$
outside some exceptional set of finite measure,
where 
\[
U(r)=\max\left\{ T(r,F),T(r,g)\right\}.
\]
By Lemma~\ref{bakerT}  and our choice of $q$
we have $T(r,g)< T(r,F)$ for large $r$.
We deduce that 
$T(r,b_m)=O(T(r,F))$ as $r\to\infty$ outside some
exceptional set of finite measure.
It now follows from Lemma~\ref{steinm} that there exists
polynomials $p_m$ such that 
\[
\sum_{m\in J}p_m M_m[g]=0.
\]
Thus $g$ satisfies an algebraic differential equation.\qed


\begin{thebibliography}{99}
\bibitem{Bak58}
I.~N.~Baker,
Zusammensetzungen ganzer Funktionen.
Math.~Z.\ {\bf 69} (1959), 121--163.

\bibitem{Bar}
D.~Bargmann,
Iteration holomorpher Funktionen.
Dissertation, Christian-Al\-brechts-Uni\-versit\"at zu Kiel, 1996.

\bibitem{Bro}
{\rm W.\ D.\ Brownawell},
{\rm On the factorization of partial differential equations}.
{\rm Can.\ J.\ Math.} {\bf 39} (1987), 825--834.


\bibitem{Clu70} J.\ Clunie, The composition of entire and meromorphic
functions. 
In  ``{\rm Mathematical essays dedicated to A.\ J.\ Macintyre},''
edited by H.\ Shankar,
Athens, Ohio:
Ohio University Press 1970, pp.~75--92.

\bibitem{Ere}
A.~E.~Eremenko,
Some functional equations connected with the iteration of 
rational functions (Russian). 
Algebra i Analiz {\bf 1} (1989), 102--116; translation in Leningrad Math.\ J.\ {\bf 1} (1990), 905--919.

\bibitem{Fat23}
{\rm P.\ Fatou},
Sur l'it\'eration analytique et les substitutions permutables.
{\rm J.\ Math.} (9) {\bf 2} (1923), 343--384.

\bibitem{Gro89}
{\rm F.\ Gross}
and
{\rm C.\ F.\ Osgood},
{\rm A simpler proof of a theorem of Steinmetz}.
{\rm J.\ Math.\ Anal.\ Appl.} {\bf 143} (1989), 490--496.
\bibitem{Gro91}
\bysame,
{\rm An extension of a theorem of Steinmetz}.
{\rm J.\ Math.\ Anal.\ Appl.} {\bf 156} (1991), 290--294.
\bibitem{Gro92}
\bysame,
{\rm Finding all solutions related to Steinmetz's theorem}.
{\rm J.\ Math.\ Anal.\ Appl.} {\bf 164} (1992), 417--421.

\bibitem{Ost}
A.\ Ostrowski,
\"Uber Dirichletsche Reihen und algebraische Differentialgleichungen.
Math.~Z.\ {\bf 8} (1920), 241--298.

\bibitem{Hay64}
{\rm W.\ K.\ Hayman},
{\rm Meromorphic functions}.
Oxford: Clarendon Press 1964.

\bibitem{Kob}
T.~Kobayashi,
Permutability and unique factorizability of certain entire functions. 
Kodai Math.\ J.\ {\bf 3} (1980), 8--25.


\bibitem{Jul22}
{\rm G.\ Julia}, M\'{e}moire sur la permutabilit\'{e}
des fractions rationnelles. {\rm Ann.\ Sci.\ \'{E}cole Norm.\
Sup.} (3) {\bf 39} (1922), 131--215.


\bibitem{Ng}
T.~W.~Ng,
Permutable entire functions and their Julia sets. 
Math.\ Proc.\ Cambridge Philos.\ Soc.\ {\bf 131} (2001), 129--138.

\bibitem{Pol26}
G.\ P\'olya,
On an integral function of an integral function.
{\rm J.\ London Math.\ Soc.} {\bf 1} (1926), 12--15.

\bibitem{Rit}
J.~F.~Ritt,
{\rm Transcendental transcendency of certain functions of Poincar\'e}.
{\rm Math.\ Ann.} {\bf 95} (1925/26), 671--682.


\bibitem{Ste80}
{\rm N.\ Steinmetz},
{\rm \"Uber faktorisierbare L\"osungen gew\"ohnlicher
Differential\-glei\-chun\-gen}.
{\rm Math.\ Z.} {\bf 170} (1980), 169--180.

\bibitem{UY}
C.-C.\  Yang and H.~Urabe,
On permutability of certain entire functions. 
J.\ London Math.\ Soc.\ (2) {\bf 14} (1976), 153--159.

\bibitem{Xio}
W.\ Xiong,
On the permutability of entire functions. 
Nanjing Daxue Xuebao Shuxue Bannian Kan {\bf 17} (2000), 56--63.

\bibitem{Zhe}
J.~H.~Zheng,
On permutability of periodic entire functions. 
J.\ Math.\ Anal.\ Appl.\ {\bf 140} (1989), 262--269.
\bibitem{ZY92}
J.~H.~Zheng and C.-C.~Yang,
Permutability of entire functions. 
Kodai Math.\ J.\ {\bf 15} (1992), 230--235.
\bibitem{ZY}
\bysame,
On the permutability of entire functions. 
J.\ Math.\ Anal.\ Appl.\ {\bf 167} (1992), 152--159.

\bibitem{ZZ}
J.~H.~Zheng and Z.~Z.~Zhou,
Permutability of entire functions satisfying certain differential equations. 
Tohoku Math.\ J.\ (2) {\bf 40} (1988), 323--330.
\end{thebibliography}
\end{document}